\documentclass[11pt]{amsart}

\newtheorem{theorem}{Theorem}[section]
\newtheorem{lemma}[theorem]{Lemma}

\theoremstyle{definition}
\newtheorem{definition}[theorem]{Definition}

\newtheorem{cor}[theorem]{Corollary}
\newtheorem{notation}[theorem]{Notation}

\newtheorem{pro}[theorem]{Proposition}

\theoremstyle{remark}

\numberwithin{equation}{section}

\begin{document}

\title{ Multidimensional continued fraction and Rational Approximation
 }

\author{{Zongduo Dai}}
\email{yangdai@public.bta.net.cn}
\thanks{This work is partly supported by NSFC (Grant No. 60173016),
 and the National 973 Project (Grant No. 1999035804)}

\author{{Kunpeng Wang}}
\email{kunpengwang@263.net}

\author{{Dingfeng Ye}}
\email{ydf@is.ac.cn}

\address{State Key Laboratory of Information Security (Graduate
School of Chinese Academy of Science), Beijing, 100039}


\date{2003,11,24,post-submission}



\begin{abstract}
 The classical continued fraction is generalized for  studying the
rational approximation problem on multi-formal Laurent series in
this paper, the construction is called  $m$-continued fraction. It
is proved that
  the approximants of an $m$-continued fraction converge to a multi-formal
Laurent series, and are best rational approximations to it;
 conversely for any multi-formal
Laurent series an algorithm called $m$-CF transform is introduced
to obtain its  $m$-continued fraction expansions; moreover, strict
$m$-continued fractions,
 which are $m$-continued fractions imposed with some additional
conditions, and  multi-formal Laurent series
    are in 1-1 correspondence.  It is shown that $m$-continued fractions can be used to study
 the multi-sequence synthesis problem.
\end{abstract}
\maketitle

\section{Introduction}
Continued fraction~\cite{Las,Rose,Stark,William,Wol}
 is a useful tool in dealing with many number
theoretic problems and in numerical computing. It is well-known
that the simple continued fraction expansion of a single real
number gives the best solution to its rational approximation
problem.
 Many people have
contrived to construct multidimensional continued fractions in
dealing with the rational approximation problem for multi-reals.
The archetypal example of a multidimensional generalization of
simple continued fraction is the Jacobi-Perron algorithm (JPA),
see \cite{Bernstein}. This algorithm and its modification are
extensively studied \cite{ito1,ito,mee,Pod}. These algorithms are
borrowed to study the same problem for multi-formal Laurent series
\cite{fengk,inoue}. But none of these algorithms guarantee best
rational approximation in general. In this paper, we deal with the
multi-rational approximation problem over the formal Laurent
series field $F((z^{-1}))$: given an element $\underline r \in
F((z^{-1}))^m$, find $\underline p \in F[z]^m, q\in F[z]$, such
that $\underline {p}/q$ approximates $\underline r$ as close as
possible while $deg(q)$ is bounded. In this setting, we give a
natural generalization of the simple 1-dimensional continued
fraction to the multidimensional case, and demonstrate that
approximants of the continued fraction expansions of a
multi-formal Laurent series $\underline r$ we define indeed give
the best rational estimates of $\underline r$.

The rational approximation problem of the multi-formal Laurent
series is motivated by sequence synthesis problem, which has
applications in the field of communication and cryptography. It is
known that the single-sequence synthesis problem can be solved by
the iterative  Berlekamp-Massey (BM) algorithm \cite{Ber,Mas}, the
Mills algorithm using continued fractions \cite{Mills}, and the
Euclidean algorithm as presented by Sugiyama Y. {\it et al.}
\cite{Sugi}.  Some consequential work shows that the continued
fraction algorithm is a powerful tool \cite{Cheng,Ni,Reed,Welch}.
Though for solving sequence synthesis problem the continued
fraction algorithm is equivalent to the BM algorithm and Euclidean
algorithm, some mysterious data structures of the latter have
natural interpretation in view of the former \cite{Dai1,Dai2}. The
BM algorithm and the Euclidean algorithm can be generalized to
solve the multi-sequence synthesis problem~\cite{Feng1,Feng2},
which can also be solved by the lattice basis reduction algorithm
as presented by Wang L.P. {\it et al.}~\cite{W,{WZ},{WZP}}. It is
natural to think of the generalization of continued fractions to
multidimensional case.

Our main contributions in this paper are: 1, we give a definition
of continued fraction in the multidimensional case; 2, set up the
transformation from continued fractions to elements of
$F((z^{-1}))^m;$ 3, prove that the approximants of a continued
fraction are best rational approximants of its corresponding
element; 4, give the procedures of computing the continued
fraction expansions of a given element of $F((z^{-1}))^m$; 5, give
a natural 1-1 correspondence between strict $m$-continued
fractions, which are $m$-continued fractions imposed with some
additional conditions, and elements of $F((z^{-1}))^m$, 6. apply
the theory to solve the multi-sequence synthesis problem.

The rest of this paper is arranged as follows: Section 2 lists the
notations used in the paper. Section 3 deals with the indexed
valuation of $F(({z^{-1}}))^m.$ Section 4 is the detailed
definition of the problem of multidimensional rational
approximation. Section 5 defines m-pre-continued fractions and
 m-continued fractions and gives some of their basic properties.
Section 6 states the main results of this paper. Section 7
completes all proofs.

\section{Notations} In this
paper we always keep the following notations.

 $F$ denotes a field,  $ F((z^{-1})) = \left\{
 \sum_{t \ge i}a_tz^{-t} \ | \ i \in Z, a_t\in F \
 \right\} $
 is the formal Laurent series field~\cite{he}, where $Z$ denotes the integer ring. $F(z)$ denotes the rational
  fraction field over $F$, $F[z]$ denotes the polynomial ring over $F$.
  The discrete valuation on $ F((z^{-1}))
 $ is denoted by $v(\cdot )$,
 {\it i.e.,} $v(\sum_{t \ge i}a_tz^{-t})=i$ if $a_i\ne 0$.
 For $\alpha= \sum_{ d\ge -t \ge 0 } a_t z^{-t} +
 \sum_{1\le t < \infty}a_t z^{-t}
\in F_2((z^{-1}))$, the polynomial part $\sum_{ d\ge -t \ge 0 }
a_t z^{-t}$ is denoted by $\lfloor \alpha \rfloor$; and the
remaining part $\sum_{1\le t < \infty}a_t z^{-t}$ is
 denoted by $\{ \alpha \}$.

$m$ denotes a positive integer,  $Z_m$ denotes the finite set
$\{1,2,\cdots,m \}$. $ F((z^{-1}))^m $ denotes the set of all
$m$-tuples
 over $ F((z^{-1}))$, similar for $F(z)^m$ and $F[z]^m$, and we always
 write those $m$-tuples in column form. For an element
$\underline{r}= (r_1(z), r_2(z),\cdots,r_m(z))^\tau$ in $
F((z^{-1}))^m$, where $\tau$ means transpose, we denote $( \lfloor
r_1(z)\rfloor, \lfloor r_2(z)\rfloor, \cdots,\lfloor
r_m(z)\rfloor)^\tau$ by $\lfloor \underline{r} \rfloor $, and \\ $
( \{ r_1(z)\}, \{ r_2(z)\},\cdots,\{ r_m(z)\})^\tau$ by $ \{
\underline{r} \}$.

 $I_n$ denotes the identity matrix of order $n$ for any positive
integer $n$, and $\underline{e}_j $ denotes the $j$-th column of
the matrix $I_n$ for $j: \; 1\le j \le n$. Denote by
$M_{a,b}(F((z^{-1})))  $ the set of all possible matrices of order
$a\times b$ over $F((z^{-1}))$ for any integers $a$ and $b$, and
similar for $M_{a,b}(F[z] )  $.

\section{Indexed valuation on $ F((z^{-1}))^m$}
 In this section the indexed valuation on $F((z^{-1}))^m$ is introduced, which will be in use
  in studying the rational approximation problem for multi-Laurent
  series.

 \begin{definition}
( { \it Linear ordering on $Z_m\times Z$ } )
 \\
For any two elements $(h,v)$ and $(h',v')$ in $Z_m \times Z$, we
say $(h,v) <(h',v')$ if $ v<v'$, or $v=v'$ and $h<h'$.
\end{definition}
It is clear that the ordering defined above on $Z_m\times Z$
 is linear \cite{Hun}, {\it i.e.}, any two elements can be compared.

\begin{lemma}\label{le:jx} Let $ (i,x)$ and $(j ,y) $ be any two elements in $ Z_m\times Z
$, then $ (i,x)<(j ,y)$ if and only if $ x+L_{i,j} \le y$, or
equivalently if and only if $ x-l_{j,i} < y $, where $ l_{i,j}=1$
 if $ i>j $, $ l_{i,j}=0$ if $ i\le j $; $L_{i,j}=1$ if $ i\ge j $, $L_{i,j}=0$ if $ i< j $.
\end{lemma}
\begin{proof}
 It is easy to check. \end{proof}

 For any $(j,x)\in Z_m\times Z$, we call the non-zero element $z^x
\underline{e}_j$ in $F((z^{-1}))^m$ a {\it monomial}, and denote
\begin{equation}\label{eq:ivm} Iv(z^x \underline{e}_j)=(j,-x)\in Z_m\times Z.
\end{equation}

\begin{notation}\label{no:M}
Denote by ${\mathcal M} $ the set of all monomials in
$F((z^{-1}))^m$. Any element $ \underline{r}$ in $F((z^{-1}))^m$
can be expressed uniquely as a
 sum of the following form: \begin{equation}\label{eq:mon}
 \underline{r}=\sum_{\underline{m}\in {\mathcal M}} c_{\underline{m}}\; \underline{m} , \;c_{\underline{m}} \in F,
\; \underline{m}\in {\mathcal M}.\end{equation} We will
 say $ \underline{m}$ is a monomial of $ \underline{r}$ and denote $ \underline{m}\in
 \underline{r}$ if $ c_{\underline{m}}\ne 0$.
  The total number ( finite or infinite ) of
the nonzero coefficients $ c_{\underline{m}} $
 will be called the {\it Hamming weight} of $\underline{r}$, denoted as
$w_H(\underline{r})=\tau$.
\end{notation}

 \begin{definition}\label{de:Iv} ({\it Indexed valuation})\\
 For an nonzero element $ \underline{r}=(r_1(z),r_2(z),\cdots,r_m(z))^\tau $ in $ F((z^{-1}))^m$,
 we call $$Iv( \underline{r} ) =\min \{ Iv(\underline{m}) \; | \; \underline{m} \in
 {\mathcal M},\; \underline{m} \in \underline{r} \; \}\in Z_m\times Z$$
the {\it indexed valuation } of $ \underline{r}$, where
$Iv(\underline{m})$ is defined as in (\ref{eq:ivm}). If $Iv(
\underline{r} ) = (h,v)$, then $v$ will be called the {\it
valuation } of $ \underline{r}$ and denoted by $v(\underline{r})$,
and $h$ will be called the {\it index } of $ \underline{r}$ and
denoted by $I(\underline{r})$. It is clear that
$$\left\{\begin{array}{l}v( \underline{r}
 ) =\min \{v( r_i(z)), 1\le i\le m \},\\
 I( \underline{r} ) =\min \{i | 1\le i\le m, v( r_i(z))=v( \underline{r} ) \}, \end{array}\right.$$
where $v( r_i(z))$ is the discrete valuation on
$F((z^{-1}))$~\cite{he}. By convention $v(0)=\infty$,
$Iv(\underline{0})=(1,\infty)$, and
$Iv(\underline{r})<Iv(\underline{0})\; \forall \underline{r}\ne
\underline{0}$.
\end{definition}

For $ \underline{r} \in F((z^{-1}))^m$ such that $w_H(
\underline{r} )= \tau < \infty$, define $$ Supp( \underline{r})=
\{ Iv(\underline{m} ) | \underline{m}\in {\mathcal M},
\underline{m}\in \underline{r} \},$$ which will be called the {\it
support} of $ \underline{r}$; and
 denote by $Supp^+(\underline{r})$ the
largest element in $Supp(\underline{r})$ if $w_H(\underline{r})<
\infty$. The following lemma is clear.
\begin{lemma}\label{le:Iv} The indexed valuation is a surjective
mapping $Iv$: $ F((z^{-1}))^m $ $\rightarrow
 (Z_m\times Z) \cup (1, \infty)$ which satisfies
 \begin{enumerate}
 \item
 $Iv(\underline{r})=(1, \infty)\Leftrightarrow \underline{r}=\underline{0}.$
 \item
 $Iv(\underline{r} u(z))=(h,v+v(u(z)) )$ for all $\underline{r}\in
 F((z^{-1}))^m$ and $u(z)\in F((z^{-1}))$, where $(h,v)=Iv(\underline{r})$.
\item
$ Iv( \underline{r}+\underline{s} ) \ge \min \{ Iv( \underline{r}
) , Iv( \underline{s}) \}$ for all $ \underline{r},
\underline{s}\in F((z^{-1}))^m $. The equality holds true if $ Iv(
\underline{r} )\ne Iv( \underline{s} ) $.
\end{enumerate}
\end{lemma}

 \section{ Rational Approximation }
 In this section we give a detailed  definition for the rational
 approximation problem on multi-formal Laurent series.
Let $ \underline{r}\in F((z^{-1}))^m$, $0\ne q(z)\in F[z]$, $
\underline{p}(z)\in F[z]^m$, we call $ \frac{
\underline{p}(z)}{q(z)}$
 a {\it rational approximant } to $ \underline{r}$
 if $v(\underline{r}-\frac{\underline{p}(z)}{q(z)} )> \deg(q(z))$, and call
 $Iv(\underline{r}-\frac{\underline{p}(z)}{q(z)})$ the {\it precision} of
 $\frac{\underline{p}(z)}{q(z)}$ to $\underline{r}$, here it is worth of pointing out that the
 components of $\frac{\underline{p}(z)}{q(z)}$
 have a common denominator $q(z)$. It is easy to see that
 $\frac{\underline{p}(z)}{q(z)}$ is a rational approximant to $ \underline{r}$
 if and only if $ \underline{p}(z)=\lfloor \underline{r}q(z) \rfloor$. So any
nonzero polynomial $q(z)$ can be the denominator of a rational
approximant to $ \underline{r}$, and such a rational approximant
is unique. As a matter of convenience, any nonzero polynomial
$q(z)$ will be called a {\it denominator } of $ \underline{r}$
with precision $Iv(\underline{r}-\frac{\lfloor\underline{r}(z)q(z)
\rfloor }{q(z)})$.

\begin{definition} ({\it Best Rational approximant}) \\
Let $ \underline{r} \in F((z^{-1}))^m$ . We call
 $\frac{\underline {p}(z)}{q(z)} \in F(z)^m$
a {\it best rational approximant } of $ \underline{r}$ if
denominators of $\underline{r}$ with degree lower than $deg(q(z))$
have precision lower than
$Iv(\underline{r}-\frac{\underline{p}(z)}{q(z)} )$, and
denominators of $\underline{r}$ with degree same as $q(z)$ have
precision no greater than
$Iv(\underline{r}-\frac{\underline{p}(z)}{q(z)} )$. In this case,
$q(z)$ is also called a {\it best denominator} of $\underline{r}$.
\end{definition}

 Let ${\mathbf S}$ be the set of all elements $ \underline{r}$
of $F((z^{-1}))^m$ with $v(\underline{r})>0$. For the rational
approximation problem on  $F((z^{-1}))^m$, we need only consider
elements in  ${\mathbf S}$. In fact, if we write $ \underline{r} =
\lfloor \underline{r} \rfloor + \{ \underline{r} \},$ then, $ \{
\underline{r} \}\in {\mathbf S}$, and $ \lfloor \underline{r}
\rfloor + \frac{ \underline{p}(z)}{q(z)} $ is a rational
approximant to $ \underline{r} $  of precision $(h,n)$, if and
only if
 $\frac{
\underline{p}(z)}{q(z)}\in F(z)^m$ is a rational approximant to $
\{ \underline{r} \} $ of the same precision $(h,n)$; and $q(z)$ is
a best denominator of $\underline{r}$ if and only if it is a best
denominator of $\{ \underline{r}\}$.

In studying the best rational approximation problem for ${
\underline{r}\in {\mathbf S}}$, motivated by the simple continued
fractions in the case $m=1$, we will introduce the
multidimensional continued fraction in the next section.

\section{ Multidimensional continued fractions }
 In this section we introduce the  $m$-pre-continued fractions and
 the $m$-continued fractions,  and give some of their basic properties.
\subsection{ $m$-pre-continued fractions}
\begin{definition} ({\it $m$-pre-continued fraction}) Let
\begin{equation}\label{eq:cf} C=[\underline{a}_0, h_1, \underline{a}_1,h_2,
\underline{a}_2,\cdots, h_k,\underline{a}_k,\cdots ], \; 0\le k<
\omega,
\end{equation}
 where $h_k\in Z_m$ and $ \underline{a}_k \in
 F[z]^m$ for all $k:\; 0\le k< \omega$, and $ \omega$ is a positive integer or
 $\infty$, then $C$ will be called
 an multidimensional-pre-continued fraction, or simply {\it $m$-pre-continued fraction},
and $\omega$ the {\it length } of $C$.
 We always assume $\underline{a}_0(z)= \underline{0}$,
 since we will see $\underline{a}_0$ is irrelevant to our concern.
\end{definition}

\begin{notation}\label{no:EI}
 Denote by
$ E_h$  the matrix of order $(m+1)$ which comes by exchanging the
$h$-th column and the $(m+1)$-th column of the identity matrix
$I_{m+1}$:
\begin{equation}\label{eq:P}
E_h=( \underline{e}_1\;\; \underline{e}_2 \cdots
\underline{e}_{h-1}\;\; \underline{e}_{m+1}\;\;
\underline{e}_{h+1} \cdots \underline{e}_{m} \;\;
\underline{e}_{h}\; )\in M_{m+1,m+1}(F[z] )  .
\end{equation}
\end{notation}

 From $C$ defined as (\ref{eq:cf}) we define
iteratively the {\it square matrices} $B_k$ of order $(m+1)$ over
$F[z]$:
\begin{equation}B_0=I_{m+1}, \; B_k= B_{k-1} E_{h_k}A( \underline{a}_k )\in M_{m+1,m+1}(F[z] ) , \; k\ge
1, \label{eq:Bk}
\end{equation}
where
\begin{equation}\label{eq:Aa} A(
\underline{a}_k) =\left(\begin{array}{cc} I_m& \underline{a}_k
\\0 &1
\end{array}\right)\in M_{m+1,m+1}(F[z] ) ; \end{equation}
and denote by $\underline{b}_k$ the last
 column of $ B_k$, and let $\underline{p}_k= \underline{p}_k(z)\in F[z]^m
 $ and $ q_k=q_k(z)\in F[z]$ be {\it components} of $\underline{b}_k$, that is
\begin{equation}
\underline{b}_k = ( \underline{p}_k, q_k )^\tau . \label{eq:pqdf}
\end{equation}

Some properties of the sequence $\{(\underline{p}_k,q_k) \}_{k\ge
0}$ are given  in the following proposition.

\begin {pro}
 \label{le:id}
\begin{enumerate}
 \item $\gcd( q_k(z),p_{k,1}(z),p_{k,2}(z),\cdots, p_{k,m}(z)
)=1$ for $k\ge 0$, where $ p_{k,j}(z)$ is the
 $j$-th component of $\underline{p}_k$.
\item\label{it:basic1+}
For $k\ge 0$, let $P_{k-1}$ be the matrix of size $m\times m$ and
$Q_{k-1}$ the matrix of size $1\times m$ such that
\begin{equation}\label{eq:Aa} B_k=\left(\begin{array}{cc} P_{k-1} &
\underline{p}_k
\\Q_{k-1} &q_k
\end{array}\right),  \end{equation}
 and denote by
 $\underline{P}_{k-1,j} (\in F[z]^m)$ the $j$th column of $P_{k-1}$, and $Q_{k-1,j}
 (\in F[z])$ the $j$th component of $Q_{k-1}$, $1\le j \le
 m$. Then for $k\ge 1$,
$$
\left(\begin{array}{c}  \underline{p}_{k} \\
 q_{k}
 \end{array}\right)
=
\left(\begin{array}{c}  P_{k-1} \\
 Q_{k-1}
 \end{array}\right)\underline{a}_k(z)+\left(\begin{array}{c}  P_{k-2,h_k} \\
 Q_{k-2,h_k}
 \end{array}\right),$$
  or explicitly
 \begin{eqnarray*} \underline{p}_k&=&
 \underline{p}_{k-1}a_{k,h_k}(z)+ \sum_{j\ne h_k,
1\le j \le m} \underline{P}_{k-2,j} a_{k,j}(z)+
\underline{P}_{k-2,h_k} ,\\
 q_k
 &=&q_{k-1}a_{k,h_k}(z)+ \sum_{j\ne h_k,
1\le j \le m} Q_{k-2,j} a_{k,j}(z)+ Q_{k-2,h_k} .
\end{eqnarray*}
 \item \label{it:basic3}
 For $1\le k < \omega$ and $j\in Z_m$, denote by
 $l(k, j)$ the largest positive integer $k'\le k$ such that
 $h_{k'}=j$, and let $l(k,j)=0$ if no such $k'$ exists.
   Then for $k\ge 1$,
   $$(\underline{P}_{k-1,j}, Q_{k-1,j})
 =\left\{ \begin{array}{ll}
 (\underline{p}_{l(k, j)-1},q_{l(k, j)-1}) & if \; l(k, j) \ge 1,\\
(\underline{e}_j, 0) & if \; l(k, j) =0. \end{array}\right.$$
\end{enumerate}
\end{pro}
\begin{proof} It is easy to prove.
\end{proof}

\subsection{Conditions 1-4 and $m$-continued fractions }
\begin{definition}\label{de:con} ({\it Conditions} 1-4)
 For the $m$-pre-continued fraction $C$ defined as (\ref{eq:cf}), we define
 condition 1 as below:
 \begin{itemize}
\item {\it Condition 1:}
 $\deg( a_{k,h_k}(z)) \ge 1, \; \forall \; 1\le k < \omega,$ where $a_{k,h_k}(z)$ denotes the
$h_k$-component of $\underline{a}_k(z)$.\end{itemize} In the
sequel, we always assume $C$ satisfies the condition 1, and
associate it with the following {\it quantities} (for each $k:
1\le k<\omega $):
\begin{eqnarray}\label{eq:vkj} \left\{ \begin{array}{ll}
t_k = \deg( a_{k,h_k}(z)), & t_0=0,\\
v_{k,j}=\sum_{h_i=j, i\le k }t_i, & v_{0,j}=0, \\
v_{k}=v_{k,h_k}, & v_0=0,
\end{array}\right.
 \end{eqnarray}
 and the {\it diagonal matrix $D_k$ } (for each $k:
1\le k<\omega $):
\begin{eqnarray}\label{eq:Dk} D_k &=& Diag. ( z^{- v_{k,1} }, z^{-
v_{k,2} }, \cdots, z^{ -v_{k,m} } )
 \\&=&\left(\begin{array}{ccc} z^{- v_{k,1} } &&\\
&\ddots&\\ &&z^{- v_{k,m} } \end{array}\right). \nonumber
 \end{eqnarray}
We define conditions 2-4 on $C$ as below (by convention
$\infty-1=\infty$):
 \begin{itemize}
\item {\it Condition 2:}
$(h_k, v_{k-1, h_k } )< ( h_{k+1}, v_{k+1} ), \; 1\le k <
\omega-1$.
\item
{\it Condition 3:} $Iv( D _k \underline{a}_k )=(h_k, v_{k-1, h_k}
),\; 1\le k < \omega$.
\item
{\it Condition 4:} $ Supp^+( D_k \underline{a}_k )<(h_{k+1},
v_{k+1} ),\; 1\le k < \omega-1$.
\end{itemize}
\end{definition}

\begin{definition}
 An $m$-pre-continued fraction $C$ is called an {\it
$m$-continued fraction} if $C$ satisfies the conditions 1-3.
 An $m$-continued
fraction $C$ is said to be {\it strict} if it satisfies the
condition 4.
\end{definition}

It is clear that Conditions $1$, $3$ and $4$ imply Condition $2$.
These conditions can be stated in some equivalent forms as shown
in the following proposition.

 \begin{pro}
\label{pro:C123}
 Let $C$, which is defined as (\ref{eq:cf}), satisfies the condition $1$, keep all notations made
 for it, and denote
\begin{equation} \label{eq:ak}
 \underline{a}_k =(a_{k,1}(z),
 \cdots , a_{k,j} (z), \cdots, a_{k,m} (z))^\tau,\; 1\le k
 <\omega,
\end{equation}
where $ a_{k,j}(z)\in F[z]$. Moreover, we associate $C$ with
 some {\it quantities } as below: for each $k: 1\le k<\omega $,
\begin{eqnarray}\label{eq:vkj2} \left\{ \begin{array}{l}d_0=0,\\
d_k = \sum_{1\le i \le k} t_i, \\
n_k = d_{k-1}+ v_k (=d_k+v_{k-1,h_k} ),
\end{array}\right.
 \end{eqnarray}
 where the notations $t_k$ and $ v_k$ are defined
 in~(\ref{eq:vkj}). Let $l_{i,j}$ and $L_{i,j}$ be defined as
 Lemma~\ref{le:jx}. Then
 \begin{enumerate} \item \label{it:C2}
 For $k: 1\le k < \omega-1$, the following conditions are equivalent:
 \begin{enumerate} \item \label{te:cc111} $(h_k, v_{k-1, h_k})<(h_{k+1}, v_{k+1})$.
 \item\label{te:cc2}
 $(h_k, n_k)<(h_{k+1}, n_{k+1})$.
\item\label{te:cc3} $ v_{k-1,h_k}-v_{k,h_{k+1}}+ L_{
h_k,h_{k+1}} \le t_{k+1}$.
 \end{enumerate}
\item\label{it:C3} For $k: 1\le k < \omega$, the following two conditions are equivalent:
 \begin{enumerate} \item \label{te:cc31}
$Iv( D_k \underline{a}_k )=( h_k, v_{k-1,h_k})$.
\item\label{te:cc32}
$ \deg(a_{k,j}(z)) \le v_{k,j}-v_{k-1,h_k} - l_{h_k,j}, \;
 \forall 1\le j \le m, j\ne h_k. $
 \end{enumerate}
\item\label{it:C4} Assume the condition 3 holds true. For $k: 1\le k < \omega-1$,
the following two conditions are equivalent:
 \begin{enumerate} \item \label{te:cc41}
 $Supp^+( D_k \underline{a}_k )<(h_{k+1},
v_{k+1} ) $. \item \label{eq:cc42}
 For every $j: 1\le j
\le m$,
 $a_{k,j}(z)$ is of the form
$$ a_{k,j}(z)=\sum_{0\le x \le
X_{k,j}}a_{k,j, x}z^{t_{k,j}-x}, \; a_{k,j,x}\in F, $$ where
 $X_{k,j}=\min\{t_{k,j}, x_{k,j} \}$, $
t_{k,j}=v_{k,j}-v_{k-1,h_k} - l_{h_k,j}$, $x_{k,j}=
v_{k+1}-v_{k-1,h_k}-l_{h_k,j}-L_{j, h_{k+1} } $.
 \end{enumerate}
\item\label{it:C23} Assume both of the condition $2$ and $3$ hold
true. For $k: \; 1\le k <\omega$ and $j\neq h_k$,
$$\left\{ \begin{array}{ll} \deg(a_{k,j}(z)) < d_k- d_{l(k,j)-1} & \quad
if \; l(k,j)\ge 1,\\
deg(a_{k,j}(z))\le 0&\quad if \; l(k,j)=0.\end{array}\right.$$
\end{enumerate}
\end{pro}
\begin{proof} It is easy to prove.
\end{proof}

\section{ Main Results}
In this section we state the main results of this paper.
\subsection{ Approximants of $m$-continued fraction-I}
 In this subsection, we always let
 $ C$, which is defined as
(\ref{eq:cf}), be an $m$-continued fraction.
 We keep
all notations made for it in section 5.

\begin{theorem}\label{th:1} $ \deg(q_k(z))=d_{k}$ for all $k: 0\le k
< \omega$.
\end{theorem}

Based on the above theorem, each pair $( \underline{p}_k, q_k)$
provides a rational fraction
$\frac{\underline{p}_k(z)}{q_{k}(z)}$, which will be called the
$k$th {\it approximant} (or, {\it convergent}) of the
$m$-continued fraction $C$.
 The properties of these approximants are summarized in the following
theorems.

\begin{theorem}\label{th:2}
$$Iv\left( \frac{\underline{p}_{k-1}(z)}{q_{k-1}(z)}
-\frac{\underline{p}_k(z)}{q_{k}(z)} \right) =(h_k, n_k), \; 1\le
k < \omega .$$
 As a consequence, the sequence $\left\{
\frac{\underline{p}_k}{q_{k}} \right\}_{k\ge 1}$ is convergent in
the case $ \omega=\infty$. \end{theorem}

 Denote
 $$ \varphi( C)=\left\{ \begin{array}{ll}\frac{
 \underline{p}_{\omega-1}(z)}{q_{\omega-1}(z)},& if \; \omega <
 \infty,\\
\lim_{k \rightarrow\infty }\frac{ \underline{p}_{k}(z)}{q_k(z)}, &
\omega =
 \infty. \end{array}\right.
$$
We see $\varphi( C)\in   F((z^{-1}))^m$, and call $C$
 an {\it $m$-continued fraction expansion } of $\varphi(
 C)$.  Denote by ${\mathbf C}( \underline{r})$ the set of all
possible $m$-continued fraction expansions of $\underline{r}$.

\begin{cor} \label{cor:1} Let $\underline{r}=\varphi( C)$.
 For $k: 1\le k
< \omega$, we have
\begin{enumerate} \item
 $Iv
(\underline{r}-\frac{\underline{p}_{k-1}(z)}{q_{k-1}(z)}) = (h_k,
n_k)$. As a consequence,
$\underline{r}q_{k-1}(z)-\underline{p}_{k-1}(z)=\{\underline{r}q_{k-1}(z)
\}$ and $Iv (\{\underline{r}q_{k-1}(z) \} ) = (h_k, v_k)$.
\item
 $Iv
(\underline{r})=(h_1, t_1)$.
 \item $Iv(\underline{p}_k(z) )=(h_1,
-d_k+t_1)$. As a consequence, $ \deg( \underline{p}_k(z) )=
d_{k}-t_1$, where $\deg( \underline{p}_k(z))$ denotes the largest
one
 among $\deg(p_{k,j}(z))$, $1\le j \le m$.
 \end{enumerate}
 \end{cor}
\begin{proof} It is easy to prove.
\end{proof}

\begin{theorem} \label{th:3}
Let $\underline{r}=\varphi( C)$. Assume $q(z)\in F[z]$, $d_k\le
\deg(q(z))<d_{k+1}$ for some $0\le k < \omega$  $(
d_{\omega}=\infty$ if $\omega< \infty )$. Then
$$ Iv (\underline{r}-\frac{\lfloor
\underline{r}q(z) \rfloor}{q(z)} ) \le Iv
(\underline{r}-\frac{\underline{p}_{k}(z)}{q_{k}(z)}).$$ As a
consequence, each $\frac{\underline{p}_{k}(z)}{q_{k}(z)}$, $0\le k
< \omega$,
 is a best rational approximant to $\underline{r}$, and the degree of any best
 denominator of $\underline{r}$ must be $d_k$ for some $k: 0 \le k < \omega$.
 \end{theorem}

\vskip3mm \noindent
\subsection{$m$-CF Transform and $m$-continued fractions-I}
Given an element $\underline{r}\in F((z^{-1}))^m$ with
$\underline{r}\ne \underline{0}$ and $v(\underline{r})>0$. In this
subsection we introduce a transform, which we call  $m$-CF
Transform, which may produce the $m$-continued fraction expansions
of $\underline{r}$.

\begin{definition} ({\it $D$-matrix})
 We call a diagonal matrix  over $F((z^{-1}))$ a
{\it $ D$-matrix } if each of its diagonal elements is a power of
$z$.
 \end{definition}

 \vskip3mm \noindent
 {\bf $m$-CF Transform:} Given $\underline{r}\in F((z^{-1}))^m$, $\underline{r}\ne \underline{0}$
 and
 $v(\underline{r})>0$. Initially,
take $\underline{a}_0=\underline{0}\in F[z^{-1}]^m$,
$\Delta_{-1}=I_m$ (the identity matrix of order $m$), and $
\beta_0=\underline{r}$.
  Suppose  $[\underline{a}_0, h_1, \underline{a}_1,h_2,
\underline{a}_2,\cdots, h_{k-1},\underline{a}_{k-1} ]$,
  $\Delta_{k-2}=Diag.( \cdots, z^{-c_{k-1,j}},\cdots  )$ which is a $D$-matrix of order $m$, and $
\beta_{k-1}=$\\ $( \cdots, \beta_{k-1,j}, \cdots)^\tau\in
F((z^{-1}))^m$ have been defined for an integer $k\ge 1$. If
$\beta_{k-1}= \underline{0}$, the algorithm terminates
 and denote $ \omega=k$. If
$\beta_{k-1}\ne \underline{0}$, then do the following steps:
\begin{enumerate}
\item Take  $ h_k = I( \Delta_{k-2}\beta_{k-1} )\in Z_m$.
\item Take
 $ \Delta_{k-1}=Diag.( \cdots, z^{-c_{k,j}},\cdots )$,
 where $c_{k,j}=c_{k-1,j}$ if $j\ne h_k$, and $c_{k,h_k}=v(\Delta_{k-2}\beta_{k-1}
 )\in Z$.
 \item Take  $\underline{a}_k= \lfloor\rho_k \rfloor - \underline{\epsilon}_k$ and  $\beta_k= \{ \rho_k \}
 +\underline{\epsilon}_k$,
 where $\rho_k=( \cdots, \rho_{k,j}, \cdots )^\tau\in F((z^{-1}))^m$, $\rho_{k,j}\in F((z^{-1}))$,
 $\rho_{k,j}=\frac{
\beta_{k-1,j}}{\beta_{k-1,h_k}}$ if $j\ne h_k$,
$\rho_{k,h_k}=\frac{ 1}{\beta_{k-1,h_k}}$, $\beta_{k-1,j}$ is the
$j$-th component of $\beta_{k-1}$, and
  $\underline{\epsilon}_k\in F[z]^m$ is chosen freely, except that it satisfies the
  condition $Iv(\Delta_{k-1} \{ \rho_k \} )\le
Iv(\Delta_{k-1}\underline{\epsilon}_k) $.
\end{enumerate}
 Denote $ \omega=\infty$ if the above procedure never terminates.

 An $m$-CF transform on $\underline{r}$ will
result in an $m$-pre-continued fraction: $C=[\underline{0}, h_1,
\underline{a}_1,h_2, \underline{a}_2,\cdots,
h_k,\underline{a}_k,\cdots] $.
 Note that $\underline{\epsilon}_k$ may not be unique at each step $k\ge 1
 $,
different choices of $\underline{\epsilon}_k$ will give different
$C$.

\begin{definition}
We denote by ${\mathbf T}( \underline{r})$ the set of all possible
$m$-pre-continued fractions obtained from $\underline{r}$ by
$m$-CF transforms.\end{definition}

\begin{definition} ({\it $\Delta$-polynomial part of $\underline{r}$}) Let $\underline{r}\in
F((z^{-1}))^m$, $\Delta$ be a D-matrix. Write
 $A=\lfloor \underline{r} \rfloor
 = \sum_{\underline{m}\in {\mathcal M}
} c_{\underline{m}}\; \underline{m}$, $c_{\underline{m}}\in F$,
$\alpha=\{ \underline{r} \}$, where $\mathcal M$ is
 defined in Notation~\ref{no:M}. We call
$$\lfloor
\underline{r} \rfloor^-_{\Delta}=\sum\limits_{\underline{m}\in
Supp_{\Delta,\alpha,-}(A)} c_{\underline{m}}\; \underline{m} $$
 the
$\Delta$-{\it polynomial part } of $\underline{r}$, where
\begin{eqnarray} Supp_{\Delta,\alpha,-}(A)&=&\{\; \underline{m}\in
{\mathcal M} \; |\; \underline{m}\in A, \; Iv(\Delta
\underline{m})< Iv(\Delta \alpha ) \; \}. \label{eq:supp+-}
\end{eqnarray}
 \end{definition}

 Denote
 $\lfloor \underline{r} \rfloor^+_{\Delta}=
\sum\limits_{\underline{m}\in Supp_{\Delta,\alpha,+}(A)}
c_{\underline{m}}\; \underline{m}$, where
\begin{eqnarray} Supp_{\Delta,\alpha,+}(A)&=&\{\; \underline{m}\in
{\mathcal M} \; |\; \underline{m}\in A, \; Iv(\Delta \alpha
)<Iv(\Delta \underline{m}) \; \}. \label{eq:supp+-2}
\end{eqnarray}
We see $\lfloor \underline{r} \rfloor=\lfloor \underline{r}
\rfloor^-_{\Delta}+
 \lfloor \underline{r} \rfloor^+_{\Delta},$ by noting that the mapping $
\underline{m}\mapsto \Delta\underline{m}$ is a bijection on
 $\mathcal M$.

\begin{lemma}\label{le:Dpart} Let $\underline{r}\in F((z^{-1}))^m$,
 $A=\lfloor
 \underline{r}\rfloor$, $\alpha=\{\underline{r}\}$; and let $\Delta$ be a
$D$-matrix. Assume $\underline{\epsilon}\in F[z]^m$, and
$\underline{a}=A-\underline{\epsilon}$. Then the following
conditions are equivalent: \begin{enumerate}\item
 $\underline{a}=\lfloor
\underline{r} \rfloor^-_\Delta $.
\item
$Iv(\Delta\alpha)\le Iv(\Delta \underline{\epsilon})$ if $
\underline{a}= \underline{0}$, and
$Supp^{+}(\Delta\underline{a})<Iv(\Delta\alpha)\le Iv(\Delta
\underline{\epsilon})$ if $ \underline{a}\ne \underline{0}$.
\item
$\left\{ \begin{array}{l}Supp( \underline{a}) \cup Supp(
\underline{\epsilon})=Supp(A),\\
Supp( \underline{a}) \subseteq Supp_{\Delta,\alpha,-}(A),\\
Supp( \underline{\epsilon}) \subseteq
Supp_{\Delta,\alpha,+}(A),\end{array}\right.$
\end{enumerate}
where $Supp_{\Delta,\alpha,-}(A)$ and $Supp_{\Delta,\alpha,+}(A)$
are defined as (\ref{eq:supp+-}) and (\ref{eq:supp+-2}). In
particular,
 $\lfloor
\underline{r} \rfloor^-_\Delta=\underline{r}$ for any
$\underline{r}\in F[z]^m$.
\end{lemma}
\begin{proof} It is easy to prove.
\end{proof}

\begin{definition}\label{de:psi} ({\it Mapping $ \psi $})
Define $ \psi ( \underline{r})$ to be the element in ${\mathbf T}(
\underline{r} )$, obtained from $\underline{r}$ by the $m$-CF
transform, where at each step we choose $\underline{a}_k=
 \lfloor R_{k-1}^{-1} \underline{s}_{k-2} \rfloor^-_{\Delta _{k-1}},$ {\it
 i.e.}, we choose
 $ \underline{\epsilon}_k=
 \lfloor
R_{k-1}^{-1} \underline{s}_{k-2} \rfloor^+_{\Delta_{k-1}}$.
\end{definition}

 The properties of $m$-CF transform are summarized in the following
theorem.

\begin{theorem}\label{th:CFrho}
Denote
 by ${\mathbf S}$ the set of all
possible $\underline{r}\in F((z^{-1}))^m$ with
$v(\underline{r})>0$, and by ${\mathbf C}^* $ the set of all
possible strict $m$-continued fractions.
 Then
 \begin{enumerate}
 \item
 ${\mathbf T}(\underline{r})={\mathbf C}(\underline{r})\; \forall \underline{r}\in {\mathbf S}$.
 \item The
 mapping
 $ \psi$ is a bijection from ${\mathbf S} $ onto ${\mathbf C}^*$,
 and $ \varphi$ is its inverse. \end{enumerate}
 \end{theorem}

\subsection{An application to multi-sequence synthesis problem}
In this subsection we show $m$-continued fractions can be used to
solve the multi-sequence synthesis problem. The multi-sequence
synthesis problem was stated in terms  of linear feedback shift
registers~\cite{Feng1,Feng2}.
  For convenience, we restate it by means of the indexed
valuation on $F((z^{-1}))^m$.  For any given sequence $r=\{c_t
\}_{t\ge 0}$ over $F$, where $c_t\in F$, we identify it with the
formal Laurent series $r(z)=\sum_{t\ge 0}c_tz^{-1-t}$ with
valuation larger than $0$, and  let
$r^{(n)}=(c_0,c_1,\cdots,c_{n-1})$ be  the length $n$ prefix of
the sequence $r$.
 For any given multi-sequences $\underline{r}=(r_1,\cdots,r_j,\cdots,r_m)$, where
 each $r_j=\{c_{j,t} \}_{t\ge 0}$  is a sequence over $F$,  we identify it with the
element $\underline{r}=(r_1(z), r_2(z), \cdots, r_m(z))^\tau \in
F((z^{-1}))^m$ with valuation larger than $0$, where
$r_j(z)=\sum_{t\ge 0}c_{j,t}z^{-1-t}$ is the formal Laurent series
identified with the $j$th sequence $r_j$, and let
$\underline{r}^{(n)}=(r_1^{(n)}, r_2^{(n)}, \cdots,
r_m^{(n)})^\tau$ be the length $n$ prefix  of the multi-sequences
$\underline{r}$. Given  a polynomial $q(z)$ over $F$, we call it a
characteristic polynomial of $\underline{r}^{(n)}$ if
$Iv(\underline{r}-\frac{\lfloor \underline{r}q(z)  \rfloor }{q(z)}
 )>( m,n )$, or equivalently
$Iv(\{ \underline{r}q(z) \} )>( m,n-\deg(q(z)) )$; and call it a
minimal polynomial of $\underline{r}^{(n)}$  if it is a
characteristic polynomial of $\underline{r}^{(n)}$ of the smallest
degree; and call
 $\deg( q(z) )$  the linear complexity of $\underline{r}^{(n)}$, denoted by $L_n( \underline{r})$,  if $q(z)$
  is a
minimal polynomial of $\underline{r}^{(n)}$. The  multi-sequence
synthesis problem is: Given a multi-sequences  $\underline{r}\in
F((z^{-1}))^m$, find a minimal polynomial and the linear
complexity  of $\underline{r}^{(n)}$ for each $n\ge 1$.

This problem is solved by using $m$-continued fractions  as shown
in the following corollary.

\begin{cor} Given a multi-sequences  $\underline{r}\in
F((z^{-1}))^m$, let $C \in {\mathbf T}( \underline{r} )$. Write
$C$ in the form of (\ref{eq:cf}). Then
 $q_k(z)$ is a minimal polynomial of $\underline{r}^{(n)}$ and $L_n( \underline{r})=d_k$
 for $n_k\le n < n_{k+1}$, $0\le k < \omega$, where we
 let $n_w=\infty$ in the case $\omega< \infty$, and $n_0=1$.
\end{cor}
\begin{proof} It is easy to prove.
  \end{proof}

\section{Proofs of Main Results}
 In this section we complete the proofs for all results
of this paper.

\subsection{Approximants of $m$-continued fraction-II} .\\
 \noindent{\bf Proof of Theorem~\ref{th:1} }: Keep all the
notations made for $C$. We prove it by induction on $k$. For
$k=1$, we have $ q_1(z)= a_{1,h_1}(z)$ by Proposition~\ref{le:id},
hence $\deg(q_1(z))=t_1=d_1$.
 Assume we are done for $<k$,
where $k\geq 2$. From Proposition~\ref{le:id} we have
$$q_k(z)=q_{k-1}(z)a_{k,h_k}(z)+ Q_{k-2,h_k} + \sum_{j\ne h_k,
1\le j \le m} Q_{k-2,j} a_{k,j}(z).$$ The wanted result
$\deg(q_k(z))=d_k$ follows by observing the following facts:
\begin{itemize}
\item
$\deg(q_{k-1}(z) a_{k,h_k}(z))= d_{k-1}+t_k=d_k$ (induction
assumption).
\item
If $Q_{k-2,j}\ne 0$, from Proposition~\ref{le:id} we see
$l(k-1,j)\ge 1$, then
$\deg(Q_{k-2,j})=\deg(q_{l(k-1,j)-1})=d_{l(k-1,j)-1}$ ( by
induction assumption). In particular, $\deg( Q_{k-2,h_k} )< d_k$.
\item For $j\ne h_k$ and $Q_{k-2,j}\ne 0$, from the above
and Proposition~\ref{pro:C123}
 we get $\deg(Q_{k-2,j}a_{k,j}(z) )<d_{l(k-1,j)-1}+
d_k-d_{l(k,j)-1}$, then $\deg(Q_{k-2,j}a_{k,j}(z) )<d_k$, since $
l(k-1,j)=l(k,j)$ in this case. \hfill $\Box$
\end{itemize}

\begin{lemma} \label{le:yz} Let $ C$, which is defined as (\ref{eq:cf}), be
 an $m$--continued fraction, and $\underline{r}=\varphi(C)$, and keep all notations
made for it. Let
 \begin{eqnarray*} \underline{y}_{r,t}&=&
\frac{\underline{p}_r(z)}{q_r(z)}
-\frac{\underline{p}_t(z)}{q_{t}(z)}, \; \forall \; 0\le r
<t,\label{eq:y}\\
 \underline{z}_{r,t}&=&\underline{p}_r(z) q_{t}(z)-\underline{p}_t(z) q_{r}(z),
 \; \forall \; 0\le r<t.\label{eq:z}\end{eqnarray*}
 Then
\begin{enumerate}
\item
 $Iv(\underline{z}_{r,t})=(h, v)$, if and only if $Iv(
\underline{y}_{r,t})=(h,v+d_r+ d_t)$.
\item The following statements are equivalent to each other:
\begin{eqnarray*}
 Iv( \underline{y}_{t-1,t})&=&(h_t, n_t), \; if \; 0<t<k, \label{eq:nk*}\\
 Iv( \underline{z}_{t-1,t})&=&(h_t, v_{t}-d_t), \; if \; 0<t<k,
 \label{eq:nk**}\\
 Iv( \underline{y}_{r,t})&=&(h_{r+1}, n_{r+1}), \; if \;0\le r<t<k,
 \label{eq:nk1}\\
 Iv( \underline{z}_{r,t})&=&(h_{r+1}, v_{r+1}-d_t), \; if \; 0\le
 r<t<k.
 \label{eq:nk2}
\end{eqnarray*}
\end{enumerate}
\end{lemma}
\begin{proof} The item (1) can be verified easily by Theorem~\ref{th:1}. The item (2)
can be obtained by a routing proof. \end{proof}

\noindent {\bf Proof of Theorem~\ref{th:2} }: Keep all the
notations made for
 $C$.
 It is enough to prove that $
Iv(\underline{y}_{t-1,t})=(h_t, n_t)$ for $ t\ge 1$,
 which is
equivalent to
 \begin{equation} Iv(\underline{z}_{t-1,t})=(h_t, v_t-d_t), \;1\le
 t \label{eq:vz}
 \end{equation}
from Lemma~\ref{le:yz}. For $t=1$,
 $\underline{z}_{0,1}=-\underline{p}_1(z)$, then
 $Iv(\underline{z}_{0,1})=Iv(\underline{p}_1(z))=(h_1,d_1-t_1)=(h_1,v_1-d_1)$ by Corollary~\ref{cor:1}.
 Suppose we are done
for $t<k$, where $k\ge 2$, which together with Lemma~\ref{le:yz}
implies
 \begin{equation}
  \label{eq:vz+} Iv( \underline{z}_{r,t})=(h_{r+1}, v_{r+1}-d_t), \;0\le r < t
 <k.
 \end{equation}
 By Proposition~\ref{le:id}
we have
\begin{eqnarray}\label{eq:36+}
&& \underline{z}_{k-1,k} \\ &= & \underline{p}_{k-1} q_{k}-
\underline{p}_kq_{k-1} \nonumber
\\&=& \underline{p}_{k-1} ( q_{k-1}a_{k,h_k}(z)+
\sum_{j\ne h_k} Q_{k-2,j}a_{k,j}(z)+Q_{k-2,h_k}) \nonumber\\& & -
 ( \underline{p}_{k-1}a_{k,h_k}(z)+ \sum_{j\ne h_k}
 \underline{P}_{k-2,j}a_{k,j}(z) + \underline{P}_{k-2,h_k} q_{k-1} \nonumber\\
&=& -(
 \underline{P}_{k-2,h_k}q_{k-1}-\underline{p}_{k-1} Q_{k-2,h_k}) \nonumber \\&& - \sum_{j\ne h_k}
(
 \underline{P}_{k-2,j}q_{k-1}-\underline{p}_{k-1} Q_{k-2,j}) a_{k,j}(z) \nonumber\\&=&
 -\underline{Z}_{k,h_k}
 - \sum_{j\ne h_k}
 \underline{Z}_{k,j}a_{k,j}(z) , \nonumber
\end{eqnarray}
where
$$
\underline{Z}_{k,j}=
 \underline{P}_{k-2,j}q_{k-1}-\underline{p}_{k-1} Q_{k-2,j}
 \; \forall 1\le j \le m.
$$
By Proposition~\ref{le:id} we can get
$$ \underline {Z}_{k,j}=\left\{ \begin{array}{ll} \underline{z}_{l(k-1,j)-1,k-1}& \; if
\quad l(k-1,j)\ge 1, \\ \underline{e}_jq_{k-1}(z)& \;if \quad
l(k-1,j)=0.
\end{array}\right.
$$
Hence by (\ref{eq:vz+})  we get
\begin{equation}\label{eq:Zj1}Iv( \underline{Z}_{k,j} ) = (j,
v_{l(k-1,j)}- d_{k-1}).\end{equation} In particular, we have
\begin{eqnarray}
Iv( \underline{Z}_{k,h_k} )&= &(h_k, v_{l(k-1,h_k)}- d_{k-1})
=(h_k,v_k-d_k). \label {eq:Zk*}
\end{eqnarray}
For $j\ne h_k$, we have
\begin{eqnarray*} Iv( \underline{Z}_{k,j} a_{k,j}(z) )&\ge &(j, v_{l(k-1,j)}-
d_{k-1}-(v_{l(k,j)}-v_{l(k-1,h_k)}-l_{h_k,j} ) ) \nonumber \\
 &=&(j, v_{l(k-1,h_k)}-d_{k-1}+ l_{h_k,j} )
 \nonumber \\ &>&
 (h_k, v_{l(k-1,h_k)}-d_{k-1}) \nonumber \\ &=&
 (h_k, v_k-d_k),
 \end{eqnarray*}
where the first line comes from (\ref{eq:Zj1}) and
 Proposition~\ref{pro:C123} ($v_{l(k,j)}=v_{k,j}$), the second line
comes from $l(k,j)=l(k-1,j)$ since $j\ne h_k$, the last two lines
come from definitions. This together with (\ref{eq:36+}) and
(\ref{eq:Zk*}) leads to
 equation (\ref{eq:vz}) when $t=k$.
 \hfill
 $\Box$

\begin{lemma}\label{le:cz}
Let $ C$, which is defined as (\ref{eq:cf}), be
 an $m$-continued fraction, and $\underline{r}=\varphi(C)$, and keep all notations
made for it. Assume $c_i(z)\in F[z]$, $\deg(c_i(z))<t_{i+1}$,
$0\le i < \omega-1$. Then
\begin{enumerate}
\item $Iv(\{\underline{r}q_i(z)c_i(z)\})=(h_{i+1},
v_{i+1}-\deg(c_i(z))).$
\item
$Iv(\{\underline{r}q_i(z)c_i(z)\})\ne
Iv(\{\underline{r}q_j(z)c_j(z)\}) $ for all $ 0\le j\ne i <
\omega$ and $ c_i(z)c_j(z)\ne 0. $
\end{enumerate}
\end{lemma}
\begin{proof} It is easy to prove.
\end{proof}

 \noindent
 {\bf Proof of Theorem~\ref{th:3}: } Keep all the
notations made for
 $C$. This theorem is true in the case $k+1=\omega< \infty$, since
$\underline{r}= \frac{\underline{p}_{k}(z)}{q_{k}(z)}$ in this
case. We need only consider the case $k+1< \omega$. Denote
$d=\deg(q(z))$. We have $d_k\le d <d_{k+1} $ by the assumption. It
is enough to prove $Iv( \{\underline{r}q(z) \} )\le (h_{k+1},
n_{k+1}-d )$.
 We can write
 $q(z) = \sum_{0\le i\le k}c_{i}(z) q_{i}(z)$, where $\deg(c_i(z)) <
\deg(q_{i+1}(z))- \deg(q_i(z)) = t_{i+1}$ for $i\ge 0$ (note that
$q_0(z)=1)$ and $\deg(c_{k}(z))=d-d_{k}\ge 0$.
 Note that
$\{\underline{r}q(z)\}=\sum\limits_{0\le i\le
k}\{\underline{r}q_{i}(z)c_{i}(z)\},$ from Lemma~\ref{le:cz} and
Lemma~\ref{le:Iv} we have
\begin{eqnarray*}
Iv( \{\underline{r}q(z) \} )&=& \min \{ Iv(
\{\underline{r}q_i(z)c_i(z) \} ) \; | \; c_i(z)\ne 0 , 0\le i \le
k \;
\} \\
 &\le & Iv(
\{\underline{r}q_{k}(z)c_{k}(z) \} )
\\
 &= & (h_{k+1}, v_{k+1}-\deg(c_{k}(z)) )
 = (h_{k+1}, n_{k+1}-d ). \hfill \Box
\end{eqnarray*}

\subsection{$m$-CF Transform and $m$-continued fractions-II}
In this subsection, at first we restate the $m$-CF Transform in
terms of matrices, then we prove Theorem~\ref{th:CFrho}. Before
giving the matrix version of $m$-CF Transform we
 give some necessary definitions.

 \begin{definition} ({\it Base matrix and its  $D$-component })\\
Let $R$ be a matrix of order $m $
 over $F((z^{-1}))$, we call it a {\it base matrix } if each of the columns of $R$ is nonzero, and
 the index of the
 $j$-th column of $R$ is $j$ for each $1\le j \le m$. For a base matrix $R$,  the $D$-matrix $
\Delta= Diag.(z^{-c_1}, z^{-c_2}, \cdots, z^{-c_m}) $ will be
called the $D$-{\it component} of $R$, where $c_j$ denotes the
valuation of the $j$-th column of $R$.
 \end{definition}

\begin{pro}\label{pro:matrix}{(\bf Matrix version of  $m$-CF
Transform)} The $m$-CF Transform which is introduced in subsection
6.2. can be
restated in terms of matrices as follows:\\
 Given $\underline{r}\in F((z^{-1}))^m$, $\underline{r}\ne \underline{0}$
 and
 $v(\underline{r})>0$.
Initially let $R_{-1}=I_m$, $ \underline{r}_0= \underline{r}$.
Suppose $(-R_{k-2}\;\;\; \underline{r}_{k-1})\in
M_{m,m+1}(F((z^{-1})))$, $k\ge 1$, is defined, where $  R_{k-2}\in
M_{m,m}(F((z^{-1}))) $ and $\underline{r}_{k-1}\in F((z^{-1}))^m
$. If $\underline{r}_{k-1}= \underline{0}$, the algorithm
terminates
 and denote $ \omega=k$. If
$\underline{r}_{k-1}\ne \underline{0}$, then do the following
steps:
\begin{enumerate}
\item Let
 $ h_k = I( \underline{r}_{k-1} )\in Z_m$.
 \item
 Let $R_{k-1}$ be the matrix of order $m$ and $\underline{s}_{k-2}  \in F(( z^{-1}))^m$ such that
 $(-R_{k-1} \;\; \; \underline{s}_{k-2} )=(-R_{k-2}\;\;\;
\underline{r}_{k-1})E_{h_k}\in M_{m,m+1}(F((z^{-1})))$, where
$E_{h_k}$ is defined in Notation~\ref{no:EI}.
\item Let $R_{k-1}^{-1}
\underline{s}_{k-2} =A_k+\alpha_k\in F((z^{-1}))^m$, where $A_k=
\lfloor R_{k-1}^{-1} \underline{s}_{k-2} \rfloor$, and $ \alpha_k=
\{ R_{k-1}^{-1} \underline{s}_{k-2} \} $.
\item Let $\underline{a}_k$ be any one in $
F[z]^m$ such that
 $\underline{a}_k= A_k - \underline{\epsilon}_k$,
 where
$\underline{\epsilon}_k\in F[z]^m$, $ Iv(\Delta_{k-1} \alpha_k
)\le Iv(\Delta_{k-1}\underline{\epsilon}_k) $,
 and
 $\Delta_{k-1}$ is the $D$-component of $R_{k-1}$.
\item Let $\underline{r}_k \in F(( z^{-1}))^m$ such that
 $ (-R_{k-1} \;\;\; \underline{r}_k)=(-R_{k-1} \;\;\;
\underline{s}_{k-2} )A(\underline{a}_k)\in
M_{m,m+1}(F((z^{-1})))$.
\end{enumerate}
 Denote $ \omega=\infty$ if the above procedure never terminates.
\end{pro}
\begin{proof}
It follows  directly from the following two lemmas
(Lemma~\ref{le:L} and Lemma~\ref{le:ma2} ).
\end{proof}

\begin{lemma}\label{le:ma2}  Keep notations made in
  Proposition~\ref{pro:matrix}. Denote
$\beta^*_{k-1}=\alpha_{k-1}+\epsilon_{k-1}$.
 For $k\ge 1$,  we have
 \begin{enumerate}
\item
$\beta^*_{k-1}=\underline{0}$ if and only if
$\underline{r}_{k-1}=\underline{0}$.
\item
$R_{k-1}$ is a base matrix.
\item
$ \underline{r}_{k-1}=R_{k-2}\beta^*_{k-1} $.
\item
$v(\Delta_{k-1,j})=\left\{ \begin{array}{lll} v(\Delta_{k-2,j}) &
if &  j\ne h_k,\\v(\Delta_{k-2}\beta^*_{k-1} ) & if & j= h_k.
\end{array}\right.$
\item
Denote $R_{k-1}^{-1}\underline{s}_{k-1} =(\cdots ,  \rho^*_{k,j},
\cdots )^\tau \in F((z^{-1}))^m$, $\rho^*_{k,j}\in F((z^{-1}))$.
Then $\rho^*_{k,j} = \frac{ \beta^*_{k-1,j}}{\beta^*_{k-1,h_k}}$
if $j\ne h_k$, $\rho^*_{k,h_k}=\frac{ 1}{\beta^*_{k-1,h_k}}$,
where $\beta^*_{k-1,j}$ is the $j$-th component of
$\beta^*_{k-1}$.
\end{enumerate}
\end{lemma}
 \begin{proof}
 It is easy to check.
\end{proof}

 \begin{lemma}\label{le:L}
 Any base matrix
 $R$ is invertible, and $ Iv(R\underline{r}
)=Iv(\Delta\underline{r})$ for all $ \underline{r}\in
F((z^{-1}))^m$, where $\Delta$ is the $D$-component of $R$.
 \end{lemma}
 \begin{proof} It is easy to prove.
\end{proof}

We  divide  Theorem~\ref{th:CFrho} in the following two parts,
then prove them separately.

 {\bf   Part 1 of  Theorem~\ref{th:CFrho}}:
 ${\mathbf T}(
\underline{r})\subseteq
 {\mathbf C}( \underline{r})$, and $\psi(\underline{r})$ is a strict $m$-continued fraction.

 {\bf   Part 2 of  Theorem~\ref{th:CFrho}}:
 ${\mathbf C}( \underline{r})\subseteq
 {\mathbf T}( \underline{r})$,  Moreover, if $ C\in  {\mathbf C}( \underline{r})$ is strict,
 then $C=\psi(\underline{r})$.

 To start the proof, we make a notation
which will be used later.

\begin{notation} Let $R$ be a base matrix,
define \begin{equation} v(R)=( c_1, c_2, \cdots, c_j , \cdots,
c_m) ,\end{equation}
 where $c_j$ denotes the
valuation of the $j$-th column of $R$. We use the notation $
(h,v)>v(R) $ to denote the fact that $v > c_h$. Notations
$(h,v)\le v(R)$, $ v(R)> (h,v)$, etc. should be interpreted
similarly.
\end{notation}

\begin{lemma}\label{le:O} Let $\Delta$ be a D-matrix, $0\ne \underline{r}\in
F((z^{-1}))^m$. Then $Iv(\Delta\underline{r})>v(\Delta)$ if
$v(\underline{r})>0$ $( $ including $\underline{r} =\underline{0}
) $; and $Iv(\Delta\underline{r})\le v(\Delta)$ if $\underline{0}
\ne \underline{r}\in F[z]^m$.
\end{lemma}
\begin{proof} It is easy to prove. \end{proof}

\begin{lemma}\label{le:1-2-3-4+} Let $R$ be a base matrix, $\Delta$ be
the $D$-component of $R$, $\underline{s}\in F((z^{-1}))^m$,
 $A=\lfloor
R^{-1}\underline{s}\rfloor$, and $\alpha=\{R^{-1}\underline{s}\}$.
Let $\underline{\epsilon}\in F[z]^m$, and
$\underline{a}=A-\underline{\epsilon}$. Then the following
conditions are equivalent:
\begin{enumerate}
\item $Iv(\Delta\alpha)\le Iv(\Delta\underline{\epsilon})$ .
\item $\alpha=\underline{\epsilon}=\underline{0}$, or $Iv(\Delta\alpha)<Iv(\Delta\underline{\epsilon})$.
\item $Iv(-R\underline{a}+\underline{s})=Iv(\Delta\alpha)$.
\item $Iv(-R\underline{a}+\underline{s})>v(\Delta)$
\end{enumerate}\end{lemma}
\begin{proof} It is easy to prove.  \end{proof}

\begin{lemma}\label{le:1-2-3-4}
Let $R$ be a base matrix, $\Delta$ be the $D$-component of $R$,
$\underline{s}\in F((z^{-1}))^m$,
 $A=\lfloor
R^{-1}\underline{s}\rfloor$, and $\alpha=\{R^{-1}\underline{s}\}$.
Let $\underline{\epsilon}\in F[z]^m$, and
$\underline{a}=A-\underline{\epsilon}$. If
$v(R)>Iv(\underline{s})$ and $ Iv(\Delta\alpha)\le
Iv(\Delta\underline{\epsilon})$, then
\begin{enumerate}
\item
$Iv(\underline{s})= Iv(\Delta \underline{a})=Iv(\Delta A)$, and
$\underline{a} \ne \underline{0}$.
\item \label{it:711}
$Iv(\Delta\underline{a})<Iv(-R\underline{a}+\underline{s})>v(
\Delta) $.
\item $\underline{a}=\lfloor R^{-1}\underline{s} \rfloor^-_\Delta$ if and only if
$ Supp^+(\Delta \underline{a} ) <
 Iv( -R\underline{a}+\underline{s}
)>v(\Delta). $
\end{enumerate}
\end{lemma}
\begin{proof} It is easy to prove. \end{proof}

\begin{lemma}\label{le:lemma12} Let $\underline{r}\in F((z^{-1}))^m$,
$v(\underline{r})>0$, and keep the quantities obtained by acting
the $m$-CF transform on $\underline{r}$, say, $ h_k$, $ R_{k-1}$,
 $\underline{s}_{k-2} $, $\alpha_k$, $\Delta_{k-1}$,  $\epsilon_k$, $\underline{a}_k$, $
\underline{r}_k$ and $ \omega$, etc. Let
$$C=[\underline{a}_0, h_1, \underline{a}_1,h_2,
\underline{a}_2,\cdots, h_k,\underline{a}_k,\cdots ] \in {\mathbf
T }( \underline{r}).
$$
Then, for $k: 1\le k < \omega$ we have
\begin{enumerate}
\item $R_{k-1}$ is a base matrix, so the $m$-CF transform is well
defined.
\item $v(R_{k-1})>Iv(\underline{s}_{k-2} )
=Iv(\Delta_{k-1}\underline{a}_k)<Iv(\underline{r}_k)>v(\Delta_{k-1})$.
\item
$ \deg(a_{k,h_k}(z))\ge 1$, {\it i.e.}, $C$ satisfies the
condition 1.
\item Based on the item (3), we may
keep the quantities associated with $C$, which are defined as in
section 5, say, $t_k$, $D_k$, $v_{k,j}$, $v_k$, etc. Then
$$\left\{ \begin{array}{l}
Iv(\Delta_{k-1,j})=Iv(R_{k-1,j})=(j,v_{k,j}).\\
\Delta_{k-1}=D_k.\\
Iv(\underline{r}_{k-1})=(h_k,v_k).\\
Iv(\underline{s}_{k-2} )=(h_k,v_{k-1,h_k}).\end{array}\right.$$
\end{enumerate}
\end{lemma}

\begin{proof}

It is easy to prove. \end{proof}

\noindent {\bf Proof of Part 1 of  Theorem~\ref{th:CFrho}}:\\
 From Lemma~\ref{le:lemma12}, whose notations are kept,  we see $C$
satisfies the conditions 1, so, in order to prove $C$ is an
$m$-continued fraction, it is enough to prove $C$ satisfies the
conditions 2 and 3.
 From Lemma~\ref{le:lemma12} we can see:
\begin{itemize}
\item For $k: 1\le k < \omega-1$, we have
$(h_k, v_{k-1,h_k})=Iv(\underline{s}_{k-2}
)<Iv(\underline{r}_k)=(h_{k+1},v_{k+1})$ , which means condition 2
is satisfied.
\item For $k: 1\le k < \omega$, we have
$(h_k, v_{k-1,h_k})=Iv(\underline{s}_{k-2}
)=Iv(\Delta_{k-1}\underline{a}_k)=Iv(D_{k}\underline{a}_k)$ ,
which shows that condition 3 is satisfied.
\end{itemize}

In order to prove $C$ is strict if $C=\psi(\underline{r})$, it is
enough to prove $C$ satisfies the condition 4. From
$C=\psi(\underline{r})$ we have $\underline{a_k}=\lfloor
R_{k-1}^{-1}\underline{s}_{k-2} \rfloor^-_{\Delta _{k-1}}$. From
Lemma~\ref{le:lemma12} we see $R_{k-1}$ is a base matrix, and
$v(R_{k-1})>Iv(\underline{s}_{k-2} )$; from definition of $m$-CF
transform we see $ Iv( \Delta_{k-1} \{ R_{k-1}^{-1}
\underline{s}_{k-2} \} ) \le Iv( \Delta_{k-1}
\underline{\epsilon}_k )$, so we can apply Lemma ~\ref{le:1-2-3-4}
( taking $R=R_{k-1}, \underline{s}= \underline{s}_{k-2} ,
\underline{\epsilon}= \underline{\epsilon}_k$ ) to get $
Supp^+(\Delta _{k-1}\underline{a_k}) < Iv(-R_{k-1}\underline{a}_k
+\underline{s}_{k-2} )=Iv(\underline{r}_k)$. Then, by
Lemma~\ref{le:lemma12} we get
 $Supp^+(\Delta _{k-1}\underline{a_k})< Iv(
\underline{r}_k)= (h_{k+1},v_{k+1})$ for $1 \le k < \omega-1$,
{\it i.e.,} the condition 4 is satisfied.

Finally we prove $\varphi(C)=\underline{r}$. Keep all the
notations made for $C$ in section 5. By definition of $m$-CF
transform, we can see
$(-R_{k-1},\underline{r}_{k})=(-R_{k-2},\underline{r}_{k-1})E_{h_k}A(\underline{a}_{k})
=\cdots =(-I_m,\underline{r})B_k$. This leads to
$\underline{r}_{k}=-\underline{p}_{k}+\underline{r}q_k$
 (the last column of $B_k$ is $(\underline{p}_k, q_k)^\tau$), so
$\underline{r}-\underline{p}_{k}/q_k=\underline{r}_{k}/q_k$. By
Lemma~\ref{le:lemma12} we have $Iv(\underline{r}_k)=
(h_{k+1},v_{k+1})$, then
$Iv(\underline{r}-\underline{p}_{k}/q_k)=(h_{k+1},v_{k+1}+\deg
(q_k)=(h_{k+1},v_{k+1}+d_k)$. Therefore $\lim_{k\rightarrow
\infty}\underline{p}_{k}/q_k=\underline{r}$ since
$\lim_{k\rightarrow \infty}d_k=\infty$, in other words,
$\varphi(C)=\underline{r}$. \hfill $\Box$

\begin{lemma}\label{le:lemma-13}
Let $C=[\underline{a}_0, h_1, \underline{a}_1,h_2,
\underline{a}_2,\cdots, h_k,\underline{a}_k,\cdots ] $, which is
 defined as
(\ref{eq:cf}), be an $m$-continued fraction, and keep all
notations made for it, say, $v_k$, $D_k$, etc. Let
$\underline{r}=\varphi(C)$. Define $R_{k-1}\in
M_{m,m}(F((z^{-1})))$, $\underline{r}_k \in F((z^{-1}))^m$ and $
\underline{s}_{k-2} \in F((z^{-1}))^m$ iteratively as below (
$1\le k < \omega$ ):
\begin{equation}\label{40}
\left\{
\begin{array}{l} R_{-1}=I_m,\;\;\underline{r}_0=\underline{r}\\
(-R_{k-1} \; \;\; \underline{s}_{k-2} )=(-R_{k-2}\;\;\;
\underline{r}_{k-1})E_{h_k}\in M_{m,m+1}(F((z^{-1})))\\
(-R_{k-1} \;\;\; \underline{r}_k)=(-R_{k-1} \;\;\;
\underline{s}_{k-2} )A(\underline{a}_k)\in M_{m,m+1}(F((z^{-1}))).
\end{array}\right.\end{equation}
 Denote by
$R_{k-1,j}$ the $j$-th column of $R_{k-1}$, denote by
$\Delta_{k-1}$ the $D$-component of $R_{k-1}$, and denote $
\underline{\epsilon}_k=\lfloor R_{k-1}^{-1}\underline{s}_{k-2}
\rfloor -\underline{a}_k. $
 Then, for $1\le k < \omega$, we have
\begin{enumerate}
\item $Iv(\underline{r}_{k-1})=(h_k,v_k)$.
\item
$R_{k-1}$ is a base matrix, and $\Delta_{k-1}=D_k$.
\item
$Iv(\underline{r}_k)=Iv(-R_{k-1}\underline{a}_k+\underline{s}_{k-2}
)>v(\Delta_{k-1})$.
\item
$v(R_{k-1})>Iv(\underline{s}_{k-2} )$. \item
 $ Iv(\Delta_{k-1}\{R_{k-1}^{-1}\underline{s}_{k-2} \})\le Iv(\Delta_{k-1}\underline{\epsilon}_k)$.
\end{enumerate}
\end{lemma}
\begin{proof}
It is easy to prove.
\end{proof}

\noindent
{\bf  Proof of  Part 2 of  Theorem~\ref{th:CFrho} }\\
 Keep notations made in Lemma~\ref{le:lemma-13}. From Lemma~\ref{le:lemma-13} we have
 $I( \underline{r}_{k-1})=h_k$ and
 $ Iv(\Delta
_{k-1}\{R_{k-1}^{-1}\underline{s}_{k-2} \})\le
Iv(\Delta_{k-1}\underline{\epsilon}_k)$ for $1\le k< \omega$,
which show that (\ref{40}) is an $m$-CF transform, so $C\in
{\mathbf T}( \underline{r} )$.
 To prove $C= \psi( \underline{r} )$ when
 $C$ is strict, it is enough to show $\underline{a}_k = \lfloor
R_{k-1}^{-1}\underline{s}_{k-2} \rfloor^-_{\Delta _{k-1}}$ for $k:
1\le k < \omega $. In fact, if $C$ is strict, we have
$Supp^+(D_k\underline{a}_k)<(h_{k+1}, v_{k+1})$ for $1\le k <
\omega-1$, then
$$Supp^+(\Delta_{k-1}\underline{a}_k)<(h_{k+1}, v_{k+1})=Iv( \underline{r}_k
)=Iv(-R_{k-1}\underline{a}_k+\underline{s}_{k-2}
)>v(\Delta_{k-1})$$ by Lemma~\ref{le:lemma-13}, thus
 $\underline{a}_k = \lfloor
R_{k-1}^{-1}\underline{s}_{k-2} \rfloor^-_{\Delta _{k-1}}$
according to Lemma~\ref{le:1-2-3-4} ( taking $R=R_{k-1},
\underline{s}= \underline{s}_{k-2} , \underline{\epsilon}=
\underline{\epsilon}_k$ ), since all the conditions for applying
Lemma~\ref{le:1-2-3-4} are
 satisfied by Lemma~\ref{le:lemma-13}.
 \hfill $\Box$

\bibliographystyle{amsplain}

\end{document}